\documentclass{ps}
\usepackage{latexsym}
\usepackage{amsmath,amssymb}
\usepackage{paralist}

\newtheorem{Definition}{Definition}[section]

\newtheorem{Lemma}[Definition]{Lemma}
\newtheorem{Theorem}[Definition]{Theorem}
\newtheorem{Corollary}[Definition]{Corollary}
\newtheorem{Remark}[Definition]{Remark}

\numberwithin{equation}{section}

\newcommand{\E}{\mbox{I\negthinspace E}}

\begin{document}
\title{MODERATE DEVIATIONS FOR A CURIE-WEISS MODEL WITH DYNAMICAL EXTERNAL FIELD}
\author{Anselm Reichenbachs}\address{Ruhr-Universit\"at Bochum, Fakult\"at f\"ur Mathematik, NA 4/26, 44780 Bochum, Germany;\\ e-mail: anselm.reichenbachs@rub.de}
\date{\today}
\begin{abstract} In the present paper we prove moderate deviations for a Curie-Weiss model with external magnetic field generated by a dynamical system, as introduced by Dombry and Guillotin-Plantard in \cite{guillotin}. The results extend those already obtained in the case of a constant external field by Eichelsbacher and L\"{o}we in \cite{eichLoew}. The Curie-Weiss model with dynamic external field is related to the so called dynamic $\mathbb{Z}$-random walks (see \cite{guillotinBuch}). We also prove a moderate deviation result for the dynamic $\mathbb{Z}$-random walk, completing the list of limit theorems for this object.\end{abstract}
%
%
\subjclass[2010]{60 F 10, 60 K 35, 82 B 44, 82 B 41, 60 G 50}
\keywords{moderate deviations, large deviations, statistical mechanics, Curie-Weiss model, dynamic random walks, ergodic theory}
\maketitle
\section{Introduction}
There is a long tradition in considering mean-field models in statistical mechanics such as the Curie-Weiss models. They can be considered as an approximation of the Ising model. Even though these models involve some strong simplifications important physical phenomena can be observed. In \cite{ellisNewman1} and \cite{ellisNewman2}, 
Ellis and Newman proved limit theorems for a class of Curie-Weiss models. These results have been extended by Dombry and Guillotin-Plantard
in \cite{guillotin} to a Curie-Weiss model with random external field generated by a dynamical system. They proved a weak law of large numbers, 
a CLT and a large deviation principle for the mean magnetization of the model. The purpose of the present paper is to 
prove \textit{moderate deviation principles} (MDP for short) for the Curie-Weiss model with dynamical external field, extending results already obtained 
for a class of Curie-Weiss models with constant external field by Eichelsbacher and L\"owe in \cite{eichLoew}. \\From a technical point of view there is no
distinction between a MDP and a large deviation principle. However, a large deviation principle is normally established on the scale of a law of large numbers, while MDPs describe the probabilities on a scale between a law of large numbers
and some central limit theorem. But typically, the rate function in a large deviations regime will depend on the distribution of the 
underlying random variables, while a MDP inherits properties of both the central limit behaviour as well as the large deviation
principle: the speed of convergence to zero of the probability of an untypical event usually is exponential while the rate function does not depend
on the fine structure of the underlying distribution. Nevertheless, there are interesting expamples which show that this ``folklore'' does not always hold true. Namely, a ``breakdown'' of a moderate deviation principle has been proved for the overlap parameter in the Hopfield model in \cite{eichLoew2}, i.e. a MDP does not hold for the hole range of scalings of the overlap parameter.
\par
We consider the following physical model: For a fixed positive integer $d$ and a finite subset $\Lambda\subset\mathbb{Z}^d$ a ferromagnetic crystal is described by a \textit{configuration space} 
$\Omega_{\Lambda}=\Omega^{\Lambda}$, $\Omega$ called \textit{spin space}, 
and random variables $\Sigma_{i}^{\Lambda}:\Omega_{\Lambda}\longrightarrow\Omega$, $\Sigma_{i}^{\Lambda}(\sigma)=\sigma_i$.
$\Sigma_{i}^{\Lambda}$ is called the \textit{spin} at site $i$. 
We restrict ourselves to the classical Curie-Weiss model, where the spins take values in $\Omega=\{+1,-1\}$. 
The crystal is exposed to an external magnetic field, described by a dynamical system $S=(E, \mathcal{A}, \mu, T)$, 
i.e. a probability space $(E, \mathcal{A}, \mu)$, 
a measure-preserving transformation $T:E\rightarrow E$ und a measurable function $f:E\rightarrow [0,1]$. 
We denote by $\beta=T^{-1}>0$ the \textit{inverse temperature} und by $J>0$ a \textit{coupling constant}. 
For a spin configuration, i.e. a realization $(\Sigma_{i}^{\Lambda})_{i\in\Lambda}=(\sigma_i)_{i\in\Lambda}$ and $x\in E$ we define the \textit{Hamiltonian} (see \cite{guillotin}), 
which specifies the energy of the given configuration $\sigma=(\sigma_i)_{i\in\Lambda}$:
$$H_{\Lambda, x}(\sigma)=-\frac{\beta J}{2|\Lambda|}\left(\sum_{i\in\Lambda}{\sigma_{i}}\right)^2-\frac{1}{2}\sum_{i\in\Lambda}\log\left(\frac{f(T^ix)}{1-f(T^ix)}\right)\sigma_i.$$
The energy is due to the interaction of the spins and the force of the external magnetic field. 
The probability of observing the system in state $\sigma=(\sigma_i)_{i\in\Lambda}$ is specified by the Gibbs measure:
$$P_{\Lambda, x}(\sigma)\equiv P_{\Lambda, x,\beta}(\sigma)=\frac{1}{Z_{\Lambda, x}}\exp\left(-\left[H_{\Lambda, x}(\sigma)\right]\right) \footnote{In the case that $\beta$ is fixed we keep quiet about this dependency in the notation of the measure.}$$
The normalizing factor $Z_{\Lambda, x}$ is called \textit{partition function}. For each configuration $\sigma=(\sigma_i)_{i\in\Lambda}$ 
we define the \textit{total magnetization} $M_{\Lambda}=\sum_{i\in\Lambda}\Sigma_i^{\Lambda}$. 
Without loss of generality we set $d=1$, $\Lambda=\{1,\dots, n\}$ in the sequel and we write $n$ instead of $\Lambda$, as well as $\Sigma_i^{(n)}$, $P_{n, x}$ and $M_n$.
So we consider a spin model on the complete graph with $n$ edges. This model belongs to the class of \textit{mean field} models, 
i.e. the spatial interaction is the same for every pair of spins ($J$ konstant). $\left\{\log\left(\frac{f(T^ix)}{1-f(T^ix)}\right)\right\}_{i\geq 1}$ specifies a magnetic field, which is inhomogeneous in space.  
The special case $f\equiv{1\over2}$ corresponds to the Curie-Weiss-Modell with zero external field. 
Furthermore, any external field $\left\{g(T^ix)\right\}_{i\geq 1}$ can be considered choosing the function $f=\frac{e^g}{1-e^g}$.
\par
In the case $\beta=0$ (infinite temperature) $P_{n,x}$ is equal to the product measure
$$\prod_{i=1}^n\left(f(T^ix)\delta_1+(1-f(T^ix))\delta_{-1}\right).$$
Then $M_n$ is a sum of $Ber(f(T^ix))-$distributed random variables $\sigma_i$ and it defines a \textit{dynamic $\mathbb{Z}$-random walk}. For further details see \cite{guillotin} or \cite{guillotinBuch}.
\par
Let us recall the definition of a large deviation principle: \\
Let $E$ be a metric space, endowed with the Borel sigma-field $\mathcal{B}(E)$ and $(\gamma_n)_n$ be a sequence of positive reals with $\gamma_n\rightarrow\infty$ as $n\rightarrow\infty$. A sequence of probability measures $(\mu_n)_{n\in\mathbb{N}}$ on $E$ is said to satisfy a large deviation principle (LDP for short) with speed $\gamma_n$ and good rate function $I:E\rightarrow[0,+\infty]$ if
\begin{itemize}
 \item $I$ is lower-semicontinuous and has compact level sets $\Phi(s)=\{x\in E: I(x)\leq s\}$, $s\in E$.
 \item For every open set $G\subset E$ it holds $$\liminf_{n\rightarrow\infty}\frac{1}{\gamma_n}\log\mu_n(G)\geq-\inf_{x\in G}I(x).$$
 \item For every closed set $A\subset E$ it holds $$\limsup_{n\rightarrow\infty}\frac{1}{\gamma_n}\log\mu_n(A)\leq-\inf_{x\in A}I(x).$$
\end{itemize}
Similarly we say that a sequence of random variables $(Y_n)_{n\in\mathbb{N}}$ with values in $E$ obeys a large deviation principle with speed $\gamma_n$ and good rate function $I:E\rightarrow[0,+\infty]$ if the sequence of their distributions does. We will speak about a moderate deviation principle (MDP), whenever the scaling of the corresponding random variable is between that of an ordinary law of large numbers and that of a central limit theorem.\par
In \cite{guillotinBuch} the authors proved a LDP for the mean magnetization $M_n/n$ in the above defined Curie-Weiss model with dynamical external field. We briefly recall these statements and outline the main ideas of the proofs.
The authors introduced a dynamic random walk, which is defined in the following way: Consider a dynamical system $S=(E, \mathcal{A}, \mu, T)$, where $(E, \mathcal{A}, \mu)$ is a probability space and $T$ is a measure-preserving transformation defined on $E$. Let $f:E\rightarrow [0,1]$ be a measurable function. For each $x\in E$ denote by $\mathbb{P}_x$ the distribution of the time-inhomogeneous random walk
$$S_0=0,\quad S_n=\sum_{i=1}^n X_i \text{ for } n\geq 1$$
with step distribution
$$\mathbb{P}_x(X_i=z)=\begin{cases}f(T^ix),& \mbox{if } z=1\\ 1-f(T^ix),& \mbox{if } z=-1\\ 0, & \mbox{otherwise.} \end{cases}$$
Using the G\"artner-Ellis Theorem (Birkoff's theorem implies the needed convergence of the logarithmic moment generating function) the authors proved the following result in \cite{guillotinBuch}:\\
For $\mu$-almost every $x\in E$, the sequence $(S_n/n)_{n\in\mathbb{N}}$ satisfies a LDP with speed $n$ and good rate function
$$\Lambda_x^*(y)=\sup_{\lambda\in\mathbb{R}}\{\left\langle\lambda, y\right\rangle-\Lambda_x(\lambda)\},$$
where
$$\Lambda_x(\lambda)=\E\left(\log\left(fe^{\lambda}+(1-f)e^{-\lambda}\right)\mid\mathcal{T}\right)(x),$$
$\mathcal{T}$ being the $\sigma$-field generated by the fixed points of the transformation $T$.\par
Under further assumptions on the dynamical system one can apply a stronger version of Birkhoff's theorem (see \cite{petersen}), which states pointwise convergence against a constant instead of $\mu$-almost sure convergence. The result reads as follows:\\
Suppose that the above defined dynamical system $S=(E, \mathcal{A}, \mu, T)$ is uniquely ergodic, with compact metric space $E$, continuous transformation $T$ and continuous function $f$. Then the above LDP holds for every point $x\in E$ with deterministic rate function
$$\Lambda(\lambda)=\int_E\log\left(f(y)e^{\lambda}+(1-f)e^{-\lambda}\right)~d\mu(y).$$
\par
The authors in \cite{guillotinBuch} also prove a functional central limit theorem for the dynamic random walk. For the purpose of completeness we also give a MDP for the dynamic random walk. To this end, we consider the centered random variables $\hat{X_i}=X_i-(2f(T^ix)-1)$ and first define precisely what a MDP is in our case of partial sums of independent random variables. We say that $\frac{1}{a_n}\hat{S_n}=\frac{1}{a_n}\sum_{i=1}^n\hat{X_i}$ obeys a MDP with rate function $I$ and speed $\frac{a_n^2}{n}\rightarrow\infty$, under the quenched measure $\mathbb{P}_x$, if $(a_n)_n$ is an increasing sequence of reals such that $\frac{a_n}{\sqrt{n}}\nearrow\infty$ , $\frac{a_n}{n}\searrow\infty$ and for all $A\in\mathcal{B}(\mathbb{R})$
$$-\inf_{t\in A^°}I(t)\leq\liminf_{n\rightarrow\infty}\frac{n}{a_n^2}\log\mathbb{P}_x\left(\frac{1}{a_n}\hat{S}_n\in A\right)\leq\limsup_{n\rightarrow\infty}\frac{n}{a_n^2}\log\mathbb{P}_x\left(\frac{1}{a_n}\hat{S}_n\in A\right)\leq-\inf_{t\in \overline{A}}I(t).$$
Our deviation results will depend on the speed of convergence in Birkoff's Ergodic Theorem. One defines for every $\alpha\in[0,1]$ the class of $\mu$-integrable functions
$$\mathcal{C}_{\alpha}(S):=\left\{h:E\rightarrow\mathbb{R}:\left|\sum_{k=1}^n\left(h(T^kx)-\int_Eh~d\mu\right)\right|=o(n^{\alpha})\quad\forall x\in E\right\}.$$
We get the following theorem.
\begin{Theorem} [Moderate deviations for the dynamic random walk under the quenched measure]
\label{thmMDPdyn}
Suppose that $f(1-f)\in\mathcal{C}_1(S)$ and $a:=\int_E4f(1-f)d\mu>0$. Then for all $x\in E\,$, $\frac{1}{a_n}\hat{S_n}$ obeys a MDP with speed $\frac{a_n^2}{n}$ rate function $I(t)=\frac{t^2}{2a}$.
\end{Theorem}
\begin{Corollary}
 If the dynamical system $S$ is uniquely ergodic with compact space $E$ and if $f$ is continuous, then the assertion of Theorem \ref{thmMDPdyn} holds.
\end{Corollary}
\begin{Remark}
 Note that a MDP under the annealed measure, i.e. under the measure $\mathbb{P}(dy)=\int_E\mathbb{P}_x(dy)~d\mu(x)$, can be obtained combining Theorems 2.1 and 2.2 in \cite{dinwoodie} with Theorem \ref{thmMDPdyn} above.
\end{Remark}
\bigskip
Now the LDP for the dynamic random walk on the integers yields a LDP for the mean magnetization $M_n/n$ via Varadhan's lemma, since its distribution is absolutely continuous with respect to the distribution of the dynamic random walk  and equal to
$$\frac{dP_{n,x}^{M_n}}{d\mathbb{P}_x^{S_n}}(y)=\frac{1}{\widehat{Z}_{n,x}}\exp\left[\frac{\beta J}{2n} y^2\right],$$
where $\widehat{Z}_{n,x}=\E_x\left\lbrace\exp\left[\frac{\beta J}{2n} (S_n)^2\right]\right\rbrace$ is a normalizing constant, i.e. the integrand is a continuous and bounded function on $[-1.1]$.  So $M_n/n$ under $P_{n,x}$ obeys a LDP with speed $n$ and good rate function
\begin{equation}
\label{ratefunction}
 I_{\beta,x}(s)=\Lambda^{*}(s)-\frac{\beta J}{2}s^2-\inf_{z\in\mathbb{R}}\left\{\Lambda^{*}(z)-\frac{\beta J}{2}z^2\right\},
\end{equation}
where $\Lambda^{*}$ denotes the Fenchel-Legendre transform from above.
\par
The authors also prove central limit theorems for the associated magnetization. Analogously to the treatment in \cite{ellisNewman1} and \cite{ellisNewman2}, the asymptotic behaviour of $M_n$ depends on the extremal points of a function $G$, which is a transformation of the rate function of the above LDP for the mean magnetization and defined by
$$G(s)=\frac{\beta J}{2}s^2-\int_E L(f(y),\beta Js)~d\mu(y).$$
Furthermore, one defines for every $n\geq 1$ the function
\begin{align}
G_n(s)&=\frac{\beta J}{2}s^2-\frac{1}{n}\log\E_x(\exp(\beta JsS_n))\\
&=\frac{\beta J}{2}s^2-\sum_{i=1}^nL(f(T^ix),\beta Js),
\end{align}
where
\begin{equation}
\label{L}
L(\phi,s):=\begin{cases}\left[0,1\right]\times\mathbb{R}\rightarrow\mathbb{R}\\ (\phi, s)\mapsto\log(\phi e^s+(1-\phi)e^{-s}).\end{cases}
\end{equation}
The function $G$ is real analytic, and the set where $G$ achieves its minimum is non-empty and finite (see Th. 3.1 in \cite{guillotin}). So we denote by $g=\min\{G(s): s\in\mathbb{R}\}$ the value of the global minimum (which is nonpositive since $G(0)=0$) and by $m_1,\dots,m_r$ the minimizers of $G$. Furthermore, one defines the \textit{type} $2k_i$ and the strength $\lambda_i>0$ of the minimum $m_i$ by
\begin{align*}
2k_i&=\min\{j\geq 0\mid G^{(j)}(m_i)\neq0\}\\
\lambda_i&= G^{(2k_i)}(m_i).
\end{align*} 
Usually, multiple minima occur for values of $\beta$ larger than some critical value $\beta_c$ and this phenomenon is called a ``phase transition''. For an explicit class of dynamical systems, the authors in \cite{guillotin} can compute a critical temperature $\beta_c$ for the model. But the situation for $\beta>\beta_c$, i.e. multiple minima of $G$, seems not to be well understood as we try to outline in the following:
In \cite{ellisNewman2} Ellis and Newman proved a law of large numbers (for the Curie-Weiss model with constant external field), i.e. they showed that the mean magnetization converges weakly to a linear combination of the minima of maximal type of $G$, whose weights can be written explicitly in terms of the types and strenghts of the corresponding minima. On the other hand, an LDP for the mean magnetization of this model also yields weak convergence to the zeros of the respective rate function. Interestingly enough, Ellis et al. recently proved in \cite{bge2} (Th. A.1) by means of convex analysis, that the set of global minimizers of $G$ coincides with the set of zeros of the LDP rate function. This general theorem can also be applied to the Curie-Weiss model with dynamical external field and yields coincidence of the set of zeros of the above LDP rate function (\ref{ratefunction}) and the set of minimizers of $G$. Astonishingly, in the treatment of the Curie-Weiss model with dynamical external field in \cite{guillotin}, the authors claimed that the mean magnetization does not converge in distribution in the case of multiple minima of $G$. Nevertheless, they proved exponential equivalence of $M_n/n$ to a linear combination of the minimizers of $G$, i.e. for every continuous bounded function $h$, the expectation of $h(M_n/n)$ under $P_{n,x}$ is equivalent, as $n$ goes to infinity, to
$$\frac{\sum_{i=1}^{r}b_{i,n}h(m_i)}{\sum_{i=1}^{r}b_{i,n}}.$$
For details on the $n$-dependent weights $b_{i,n}$ see Th. 3.2 in \cite{guillotin}.

\par
For the case of a unique minimum $m$ of $G$, the following limit theorem for the fluctuations of $M_n/n$ around $m$ has been proved in \cite{guillotin} (Th. 3.3): Assume that the unique minimum $m$ of $G$ is of type $2k$ and strenght $\lambda$ and that for every $j\in\{1,\dots,2k\}$, the function $\frac{\partial^j}{\partial s^j}L(f(.),\beta Jm)$ belongs to the set $\mathcal{C}_{j/2k}(S)$. Then, the following convergence of measures holds:
$$\frac{M_n-nm}{n^{1-1/2k}}\Rightarrow Z(2k,\tilde{\lambda}),$$
where $Z(2k,\tilde{\lambda})$ is the probability measure with density function
$$C\exp\left(-\tilde{\lambda}s^{2k}/(2k)!\right),$$
$C$ being a normalizing constant and
$$\tilde{\lambda}=\begin{cases}
\left(\frac{1}{\lambda}-\frac{1}{\beta J}\right)^{-1} &, \text{ if } k=1\\ \lambda &, \text{ if } k\geq 2.
\end{cases}$$
\par
The purpose of the present paper is to analyze the asymptotic behaviour of $M_n$ on a moderate deviation scale. Our results read as follows:
\begin{Theorem}[Moderate deviations for the Curie-Weiss model with dynamical external field, conditioned version]
\label{th1}
Let $m$ be a (local or global) minimum of $G$ of type $2k$ und strength $\lambda$ and assume that for every $j\in\left\{1,\dots,2k\right\}$ the function $y\mapsto\frac{\partial^j}{\partial s^j}L(f(y),\beta Jm)$ belongs to the class $\mathcal{C}_{\frac{j}{2k}}(S)$. Then there exists an $A=A(m)>0$ such that for all $0<a<A$ and for every $1-\frac{1}{2k}<\alpha<1$ the sequence of measures
$$\left\{P_{n,x}\left(\frac{M_n-nm}{n^{\alpha}}\in\bullet\;\bigg\vert\frac{M_n}{n}\in[m-a,m+a]\right)\right\}_{n\in\mathbb{N}}$$
satisfies a MDP with speed $n^{1-2k+2k\alpha}$ and rate function
$$I(z)\equiv I_{k,\lambda,\beta,J}(z):=\begin{cases}\frac{z^2}{2\sigma^2} &,\quad k=1\\ \lambda\frac{z^{2k}}{(2k)!} &,\quad k\geq 2,\end{cases}$$ where $\sigma^2=\frac{1}{\lambda}-\frac{1}{\beta J}$.
\end{Theorem}
\begin{Theorem}[Moderate deviations for the Curie-Weiss model with dynamical external field, unconditioned version]
\label{th1a}
 Assume that $G$ has a unique global minimum $m$ of type $2k$ und strength $\lambda$ and that for every $j\in\left\{1,\dots,2k\right\}$ the function $y\mapsto\frac{\partial^j}{\partial s^j}L(f(y),\beta Jm)$ belongs to the class $\mathcal{C}_{\frac{j}{2k}}(S)$. Then for every $1-\frac{1}{2k}<\alpha<1$ the sequence of measures
$$\left\{P_{n,x}\left(\frac{M_n-nm}{n^{\alpha}}\in\bullet\right)\right\}_{n\in\mathbb{N}}$$
satisfies a MDP with speed $n^{1-2k+2k\alpha}$ and rate function
$$I(z)\equiv I_{k,\lambda,\beta,J}(z):=\begin{cases}\frac{z^2}{2\sigma^2} &,\quad k=1\\ \lambda\frac{z^{2k}}{(2k)!} &,\quad k\geq 2,\end{cases}$$ where $\sigma^2=\frac{1}{\lambda}-\frac{1}{\beta J}$.
\end{Theorem}
\section{Auxiliary results}
In this section we state several lemmas that we will need in the proofs of our main theorems. The first lemma contains some important information about the sequence of functions $G_n$ and the function $G$, as defined in the introduction of the present paper. For the proof we refer to Theorem 3.1, Lemma 3.2 und Lemma 3.4 in \cite{guillotin} respectively.
\begin{Lemma}
\label{dreiLemma}
 \begin{enumerate}[(i)]
  \item The function $G$ is real analytic and the set where $G$ attains its global minimum is non-empty end finite.
  \item The sequence of functions $(G_n)_{n\geq 1}$ converges uniformly to $G$ on compacta of $\mathbb{R}$ as $n$ goes to infinity. Furthermore, the sequence of derivative functions $(G_n^{(k)})_{n\geq 1}$ converges uniformly to $G^{(k)}$ for every $k\geq 1$ on compacta of $\mathbb{R}$ as $n$ goes to infinity.
  \item Let $A\subset\mathbb{R}$ be a closed subset containing no global minima of $G$. Then there exists $\epsilon>0$ such that
$$e^{ng}\int_Ae^{-nG_n(s)}~ds=\mathcal{O}(e^{-n\epsilon}).$$
\textit{(where $g$ is the value of the global minimum of $G$.)}
 \end{enumerate}
\end{Lemma}
\bigskip
The following Lemma is a key ingredient for the proof of our MDPs. It is based on the Taylor expansion of $G$ and a slight generalization of Lemma 3.3 in \cite{guillotin}.
\begin{Lemma}
\label{hauptlemma}
Let $m$ be a (local or global) minimum of $G$ of type $2k$ and strenght $\lambda$. Suppose that for every $j\in\{1,\dots,2k\}$ the function $y\mapsto\frac{\partial^j}{\partial s^j}L(f(y),\beta Jm)$ belongs to the set $\mathcal{C}_{\frac{j}{2k}}$. Let $1-\frac{1}{2k}<\alpha<1$. Then the following assertions hold: 
\begin{enumerate}[(i)]
\item For every $s\in\mathbb{R}$
\begin{equation}
\label{glmKonvDiff}
\lim_{n\rightarrow\infty}\frac{1}{n^{2k(\alpha-1)}}(G_n(m+sn^{-(1-\alpha)})-G_n(m))=\lambda\frac{s^{2k}}{(2k)!}.
\end{equation}
The convergence is uniform on compact intervals of the form $\left[-M,M\right]$.
\item
There exist $r>0$ and $N\in\mathbb{N}$ such that for all $n\geq N$ and $s\in\left[-rn^{(1-\alpha)}, rn^{(1-\alpha)}\right]$ the following upper bound is valid:
\begin{equation}
\label{abschUnten}
\frac{1}{n^{2k(\alpha-1)}}(G_n(m+sn^{-(1-\alpha)})-G_n(m))\geq\frac{\lambda}{2}\frac{s^{2k}}{(2k)!}-\sum_{j=1}^{2k-1}\left| s\right|^j.
\end{equation}
\end{enumerate}
\end{Lemma}
\par\noindent
\textit{Proof:} ad $(i)$: Let $s\in\mathbb{R}$ and $u=sn^{-(1-\alpha)}$. Taylor expansion yields
\begin{equation}
\label{taylorEntw}
G_n(m+u)-G_n(m)=\sum_{j=1}^{2k}\frac{G_n^j(m)}{j!}u^j+R_n(u),
\end{equation}
where the remainder $R_n$ can be written in the integral form
$$R_n(u)=R_{n,2k}(m+u)=\int_m^{m+u}\frac{((m+u)-t)^{2k}}{(2k)!}G_n^{(2k+1)}(t)~dt$$
or rather (substituting $t=m+\vartheta u$)
$$R_n(u)=\frac{u^{2k+1}}{(2k)!}\int_0^1(1-\vartheta)^{2k}G_n^{(2k+1)}(m+\vartheta u)~d\vartheta.$$
The $j$-th derivative of $G_n$ in $m$ ist equal to
$$G_n^{(j)}(m)=P_j(m)-\frac{(\beta J)^j}{n}\sum_{i=1}^n\frac{\partial^j}{\partial s^j}L(f(T^ix),\beta Jm),$$
where
$$P_j(m)= \begin{cases}\beta Jm &,\quad j=1\\ \beta J &,\quad j=2\\ 0 &,\quad \text{otherwise}.\end{cases}$$
Lemma \ref{dreiLemma} \textit{(ii)} states that for every $j\in\{1,2,\dots\}$
$$G_n^{(j)}(m)\underset{n\rightarrow\infty}{\longrightarrow}G^{(j)}(m)=P_j(m)-(\beta J)^j\int_E\frac{\partial^j}{\partial s^j}L(f(y),\beta Jm)~d\mu(y).$$
The assumption that $y\mapsto\frac{\partial^j}{\partial s^j}L(f(y),\beta Jm)$ belongs to the set $\mathcal{C}_{\frac{j}{2k}}$ implies that for all $j\in\{1,\dots,2k\}$
$$n\left| G_n^{(j)}(m)-G^{(j)}(m)\right|=o(n^{\frac{j}{2k}}).$$
Since $$G^{(j)}(m)=\begin{cases}0 &,\quad j=1,\dots, 2k-1\\ \lambda &,\quad j=2k\end{cases}$$
it follows for every $j\in\{1,\dots,2k-1\}$
\begin{equation}
\label{konv1}
G_n^{(j)}(m)n^{1-\frac{j}{2k}}\underset{n\rightarrow\infty}{\longrightarrow}0
\end{equation}
and for $j=2k$
\begin{equation}
\label{konv2}
G_n^{(2k)}(m)\underset{n\rightarrow\infty}{\longrightarrow}\lambda.
\end{equation}
For the remainder we get
$$\frac{1}{n^{2k(\alpha-1)}}R_n(sn^{-(1-\alpha)})=\frac{s^{2k+1}n^{-(1-\alpha)}}{(2k)!}\int_0^1(1-\vartheta)^{2k}G_n^{(2k+1)}(m+\vartheta sn^{-(1-\alpha)})~d\vartheta.$$
This implies
\begin{equation}
\label{konv3}
\frac{1}{n^{2k(\alpha-1)}}R_n(sn^{-(1-\alpha)})\underset{n\rightarrow\infty}{\longrightarrow}0,
\end{equation}
since $G_n^{(2k+1)}$ is uniformly bounded on the compact interval $\left[m-\left| s\right|,m+\left| s\right|\right]$. Using \eqref{konv1}, \eqref{konv2} and \eqref{konv3} we thus get
\begin{align}
\frac{1}{n^{2k(\alpha-1)}}(G_n(m+sn^{-(1-\alpha)})-G_n(m))&=\sum_{j=1}^{2k}\frac{G_n^{(j)}(m)}{j!}n^{(2k-j)(1-\alpha)}s^j\notag\\
&\quad +\frac{1}{n^{2k(\alpha-1)}}R_n(sn^{-(1-\alpha)})\label{gleichheit}\\
&\quad\underset{n\rightarrow\infty}{\longrightarrow}\lambda\frac{s^{2k}}{(2k)!}.
\end{align}
The first $2k-1$ summands fade away in the limit $n\rightarrow\infty$ since
for $j\in\{1,\dots,2k-1\}$ we have
\begin{align}
G_n^{(j)}(m)n^{(2k-j)(1-\alpha)}&=G_n^{(j)}(m)(n^{(1-\frac{j}{2k})})^{2k(1-\alpha)}\cdot\underbrace{\frac{\left(n^{(1-\frac{j}{2k})}\right)^{1-2k(1-\alpha)}}{\left(n^{(1-\frac{j}{2k})}\right)^{1-2k(1-\alpha)}}}_{=1}\notag\\
&=\frac{G_n^{(j)}(m)n^{1-\frac{j}{2k}}}{n^{(1-\frac{j}{2k})(1-2k(1-\alpha))}}\notag\\
&\quad\underset{n\rightarrow\infty}{\longrightarrow}0\label{konv3}.
\end{align}
here the nominator converges to $0$ (see \eqref{konv1}) and the exponent of $n$ in the denominator
is positive, so it converges to $+\infty$ (since the assumption $1-\frac 1{2k}<\alpha<1$ yields $(1-\frac{j}{2k})(1-2k(1-\alpha))>0$). We thus have proved the first assertion \eqref{glmKonvDiff} of the Lemma. The convergence is uniform for $s\in\left[-M,M\right]$, since $s$ is independent of $n$ on the right hand side of equation \eqref{gleichheit}.
\par
ad $(ii)$: We now prove \eqref{abschUnten}. Equations \eqref{konv3} and \eqref{konv2} imply, that there exists $N\in\mathbb{N}$, such that for all $n\geq N$ and all $j\in\{1,\dots,2k-1\}$
\begin{equation}
\label{absch1}
\left|\frac{G_n^{(j)}(m)n^{(2k-j)(1-\alpha)}}{j!}\right|\leq 1
\end{equation}
and for $j=2k$
$$G_n^{(2k)}(m)\geq\frac{3}{4}\lambda.$$ Now, because of the uniform convergence of $G_n^{(2k+1)}$ on compacta of $\mathbb{R}$ there exists $M>0$ such that $\left| G_n^{(2k+1)}(s)\right|\leq M$ for all $n\in\mathbb{N}$ and for all $s\in\left[m-1,m+1\right]$. Let $r:=\min\big\{\frac{(2k+1)\lambda}{4M},1\big\}$. Then for all $n\in\mathbb{N}$ and for all $s\in\left[-rn^{(1-\alpha)},rn^{(1-\alpha)}\right]$ we have
\begin{align}
-\frac{(2k)!}{s^{2k}}\frac{1}{n^{2k(\alpha-1)}}R_n(sn^{-(1-\alpha)})&\leq\left| sn^{-(1-\alpha)}\int_0^1(1-\vartheta)^{2k}G_n{(2k+1)}(m+\vartheta sn^{-(1-\alpha)})~d\vartheta\right|\notag\\
&\leq\left| sn^{-(1-\alpha)}\right|\cdot M\left[-\frac{1}{2k+1}(1-\vartheta)^{(2k+1)}\right]_0^1\notag\\
&\leq r\cdot M\cdot\frac{1}{2k+1}\notag\\
&\leq\frac{\lambda}{4}\label{absch2},
\end{align}
thus
$$\frac{1}{n^{2k(\alpha-1)}}R_n(sn^{-(1-\alpha)})\geq-\frac{\lambda}{4}\frac{s^{2k}}{(2k)!}.$$
Using the estimates \eqref{absch1} and \eqref{absch2} equation \eqref{gleichheit} yields for all $n\geq N$
\begin{align*}
&\quad\frac{1}{n^{2k(\alpha-1)}}(G_n(m+sn^{-(1-\alpha)})-G_n(m))\\
&=\sum_{j=1}^{2k}\frac{G_n^{(j)}(m)}{j!}n^{(2k-j)(1-\alpha)}s^j+\frac{1}{n^{2k(\alpha-1)}}R_n(sn^{-(1-\alpha)})\\
&\geq\sum_{j=1}^{2k-1}-\left|\frac{G_n^{(j)}(m)}{j!}n^{(2k-j)(1-\alpha)}\right|s^j+\frac{G_n^{(2k)}(m)s^{2k}}{(2k)!}-\frac{\lambda}{4}\frac{s^{2k}}{(2k)!}\\
&\geq\frac{\lambda}{2}\frac{s^{2k}}{(2k)!}-\sum_{j=1}^{2k-1}\left|s\right|^j
\end{align*}
This is assertion $(ii)$ of the Lemma. The proof is complete.\par\hspace*{\fill}$\Box$\par\vskip2ex
\bigskip
The following lemma concerns a well known transformation of our mean-field measure, sometimes called the Hubbard-Stratonovich transform in the literature. For the proof we refer to Lemma 3.1 in \cite{guillotin}.
\begin{Lemma}
\label{trafo}
 Let $W$ be a $\mathcal{N}(0, \frac{1}{\beta J})$-distributed random variable, defined on some probability space $(\Omega, \mathcal{F}, Q)$ and independent of $M_n$ for every $n\geq 1$, and let $m$ and $\alpha$ be some real numbers. Then the random variable
$$\frac{M_n-nm}{n^{\alpha}}+\frac{W}{n^{\alpha-1/2}}$$
under the measure $Q_{n,x}:=P_{n,x}\otimes Q$ has a density with respect to the Lebesgue measure given by
\begin{equation}
\label{gewuenscht}
\frac{\exp(-nG_n(m+sn^{-(1-\alpha)}))}{\int_{\mathbb{R}}\exp(-nG_n(m+sn^{-(1-\alpha)}))ds}.
\end{equation}
\end{Lemma}
The usefulness of the previous lemma lies in the fact that one can often prove MDPs for the convolution, using the Taylor series expansion of $G$. Clearly, the type $2k$ and strenght $\lambda$ of the global minimum $m$ of $G$ will therefore play an important role. We next state two lemmas which ensure that it does not matter whether we consider the sequence of measures $P_{n,x}\circ\left(\frac{M_n-nm}{n^{\alpha}}\right)^{-1}$ or the sequence $P_{n,x}\otimes Q\circ\left(\frac{M_n-nm}{n^{\alpha}}+\frac{W}{n^{\alpha-1/2}}\right)^{-1}$ as long as $k\geq 2$.
\begin{Lemma}
\label{lemmaExp1}
 If the sequence of random variables $\frac{M_n-nm}{n^\alpha}+\frac{W}{n^{\alpha-1/2}}$ satisfies a MDP with respect to $Q_{n,x}=P_{n,x}\otimes Q$ with speed $n^{\gamma}$, $\gamma<2\alpha-1$ and rate function $I$, then so does $\frac{M_n-nm}{n^\alpha}$ with respect to $P_{n,x}$ and the speed and rate function agree.
\end{Lemma}
\noindent
\textit{Proof:} The proof is based an exponential equivalence and can be found in \cite{eichLoew}. Nevertheless we give the proof in order to allude to the problems that arise in the case $k=1$. One shows that the two sequences $\frac{M_n-nm}{n^{\alpha}}+\frac{W}{n^{\alpha-1/2}}$ and $\frac{M_n-nm}{n^\alpha}$ are exponentially equivalent and therefore have the same moderate deviation behaviour (see Th. 4.2.13 in \cite{dembo}). For all  $\epsilon>0$ the following estimate holds:
\begin{align*}
&\quad P_{n,x}\otimes Q\left(\left|\frac{M_n-nm}{n^\alpha}+\frac{W}{n^{\alpha-1/2}}-\frac{M_n-nm}{n^\alpha}\right|>\epsilon\right)=P_{n,x}\otimes Q\left(\left|\frac{W}{n^{\alpha-1/2}}\right|>\epsilon\right)\\
&=Q\left(\left|W\right|>\epsilon n^{\alpha-\frac{1}{2}}\right)\leq\sqrt{\frac{2\beta J}{\pi}}\frac{1}{\epsilon n^{\alpha-\frac{1}{2}}}\exp\left(-\frac{\beta J}{2}\epsilon^2n^{2\alpha-1}\right).
\end{align*}
This implies
\begin{align*}
&\quad\limsup_{n\rightarrow\infty}\frac{1}{n^\gamma}\log P_{n,x}\otimes Q\left(\left|\frac{M_n-nm}{n^\alpha}+\frac{W}{n^{\alpha-1/2}}-\frac{M_n-nm}{n^\alpha}\right|>\epsilon\right)\\
&\leq\limsup_{n\rightarrow\infty}\left(\frac{-\frac{\beta J}{2}\epsilon^2n^{2\alpha-1}-\log\left(\sqrt{\frac{\pi}{2\beta J}}\epsilon n^{\alpha-\frac{1}{2}}\right)}{n^{\gamma}}\right)\\
&=\limsup_{n\rightarrow\infty}\left(-\frac{\beta J}{2}\epsilon^2n^{(2\alpha-1)-\gamma}-\frac{\log\left(\sqrt{\frac{\pi}{2\beta J}}n^{\alpha-\frac{1}{2}}\right)}{n^{\gamma}}\right)\\
&=-\infty,
\end{align*}
since $\gamma<2\alpha-1$ by assumption.\par\hspace*{\fill}$\Box$\par\vskip2ex
\par
In the case $k=1$ the speed of the MDP in Theorems \ref{th1} and \ref{th1a} is of the same order as the variance of the respective Gaussian random variable. Therefore the above argument for exponential equivalence fails. In \cite{eichLoew} the authors proved a ``transfer principle'' for LDP which can be applied in this special case (see Proposition A.1 in \cite{eichLoew}). The application of this Proposition reads as follows:
\begin{Lemma}
\label{lemmaExp2}
 Suppose that $\frac{M_n-nm}{n^\alpha}+\frac{W}{n^{\alpha-1/2}}$ satisfies a MDP with respect to $Q_{n,x}=P_{n,x}\otimes Q$ with speed $n^{2\alpha-1}$ and rate function $\lambda\frac{z^2}{2}$ for a $\lambda\neq\beta J$. Then so does $\frac{M_n-nm}{n^\alpha}$ with respect to $P_{n,x}$, $1-\frac{1}{2}<\alpha<1$, with the same speed and rate function $\frac{y^2}{2\sigma^2}$, where $\sigma^2=\frac{1}{\lambda}-\frac{1}{\beta J}$. 
\end{Lemma}
\noindent
\textit{Proof:} See Lemma 3.6 in \cite{eichLoew}.
\bigskip
\par
Our last lemma can be considered as a starting point of the \textit{Laplace method} in the theory of large deviations. It will be used in the proof of our main theorems. Though it is often used implicitly in the literature we could not find a proof in the relevant books. We therefore give an own proof.
\begin{Lemma}
 \label{laplaceLemma}
Let $f:\mathbb{R}\rightarrow\mathbb{R}$ be a continuous function, $M>0$ a real number and $\gamma_n\rightarrow\infty$ a sequence of positive integers. Then
$$\lim_{n\rightarrow\infty}\frac{1}{\gamma_n}\log\int_{\{\mid x\mid\leq M\}}\exp\left[\gamma_nf(x)\right]dx=\max_{\{\mid x\mid\leq M\}}f(x).$$
\end{Lemma}
\noindent
\textit{Proof:} Since $f$ is continuous it attains its supremum on the compact interval $[-M,M]$. We get the following estimates for the limes superior and the limes inferior: 
\begin{align*}
\limsup_{n\rightarrow\infty}\frac{1}{\gamma_n}\log\int_{\{\mid x\mid\leq M\}}\exp\left[\gamma_nf(x)\right]dx&\leq\limsup_{n\rightarrow\infty}\frac{1}{\gamma_n}\log\int_{\{\mid x\mid\leq M\}}\exp\left[\gamma_n\max_{\mid y\mid\leq M}f(y)\right]dx\\
&=\limsup_{n\rightarrow\infty}\frac{1}{\gamma_n}\log\left\{\exp\left[\gamma_n\max_{\mid y\mid\leq M}f(y)\right]\cdot 2M\right\}\\
&=\max_{\mid y\mid\leq M}f(y)
\end{align*}
For $\epsilon>0$ and $y\in[-M,M]$ let $O_{y, \epsilon}:=\{x\in[-M,M]\mid f(x)>f(y)-\epsilon\}$. Since $f$ is continuous this level set is measurable with respect to Lebesgue measure and for every $\epsilon>0$ there exists a $\delta>0$ such that $B(y,\delta)\cap [-M,M]$ is contained in $O_{y, \epsilon}$, so this set has positive Lebesgue measure. Now let $\epsilon>0$ and $y\in[-M,M]$. Then
\begin{align*}
\int_{\{\mid x\mid\leq M\}}\exp\left[\gamma_nf(x)\right]dx&\geq\int_{O_{y, \epsilon}}\exp\left[\gamma_nf(x)\right]dx\\
&\geq\int_{ O_{y, \epsilon}}\exp\left[\gamma_n(f(y)-\epsilon)\right]dx\\
&=\exp\left[\gamma_n(f(y)-\epsilon)\right]\cdot\underbrace{\lambda(O_{y, \epsilon})}_{>0}.
\end{align*}
This implies
$$\liminf_{n\rightarrow\infty}\frac{1}{\gamma_n}\log\int_{\{\mid x\mid\leq M\}}\exp\left[\gamma_nf(x)\right]dx\geq f(y)-\epsilon$$
for all $y\in [-M,M]$ and all $\epsilon>0$. In particular we get
$$\liminf_{n\rightarrow\infty}\frac{1}{\gamma_n}\log\int_{\{\mid x\mid\leq M\}}\exp\left[\gamma_nf(x)\right]dx\geq\max_{\mid y\mid\leq M}f(y).$$\par\hspace*{\fill}$\Box$\par\vskip2ex

\section{Proofs}
In this section we first give the proof of our moderate deviation result \ref{thmMDPdyn} for the dynamic random walk. Thereafter we prove our main Theorems \ref{th1} and \ref{th1a}, i.e. MDPs for the fluctuations of the mean magnetization around the minimizer of $G$. Our proofs use \textit{Laplace method}, an equivalent formulation of a LDP which is based on the asymptotic analysis of the scaled logarithms of certain expectations. We refer to the book \cite{dupuis} for a detailed introduction to this approach to large deviation theory. The Laplace method has been successfully applied in the context of the Blume-Emery-Griffiths model in \cite{bge}, where the authors prove MDP's for the case of size dependent temperatures, i.e. the limit results are obtained as the pair $(\beta, J)$ converges along appropriate sequences $(\beta_n, J_n)$ to points belonging to various subsets of the phase digram. Our proof has been inspired by this approach since we have a similar $n$-dependence in the Hubbard-Stratonovich transform via the functions $G_n$.\par\noindent
\textit{Proof of Theorem \ref{thmMDPdyn}:}
In order to apply G\"{a}rtner-Ellis-Theorem \cite{dembo} we consider the Laplace transform of our random variable of interest. Denote by $s_n^2=\sum_{i=1}^n\text{Var}_x(Y_i)=\sum_{i=1}^n4f(T^ix)(1-f(T^ix))$ and let $x\in E$. Let us fix $x\in E$ and write $\lambda_i\equiv\lambda_i(x)=f(T^ix)$. We then get for each $t\in\mathbb{R}$
\begin{align*}
 \log\E_x\left[e^{t\frac{a_n^2}{n}(\hat{S}_n/a_n)}\right]&=\sum_{i=1}^n\log\E_x\left[ta_n(X_i-2\lambda_i+1)/n\right]\\
&=\sum_{i=1}^n\log\left(\lambda_ie^{ta_n(X_i-2\lambda_i+1)/n}+(1-\lambda_i)e^{ta_n(-2\lambda_i)/n}\right)\\
&=\sum_{i=1}^n\log\left(e^{-ta_n2\lambda_i/n\left(\lambda_ie^(2ta_n/n+1-\lambda_i)\right)}\right)\\
&=\frac{-ta_n2\sum_{i=1}^n\lambda_i}{n}+\sum_{i=1}^n\log\left(1+\lambda_i(e^{2ta_n/n}-1)\right)
\end{align*}
For $n$ large enough we have $\lambda_i(e^{2ta_n/n}-1)\in[0,1]$ so that we can use Taylor expansion for the logarithm.
\begin{align}
 \log\left(1+\lambda_i\left(e^{2ta_n/n-1}\right)\right)&=\lambda_i(e^{2ta_n/n}-1)-\frac{\lambda_i^2}{2}(^{2ta_n/n-1})^2+\frac{\lambda_i^3}{3}(e^{2ta_n/n-1})^3+o\left(\lambda_i^4(e^{2ta_n/n-1})^4\right)\notag\\
&=\lambda_i\frac{2ta_n}{n}+\lambda_i\frac{(2t)^2a_n^2}{2n^2}+\lambda_i\frac{(2t)^3a_n^3}{3!n^3}+o\left(\frac{\lambda_ia_n^4}{n^4}\right)-\lambda_i^2\frac{(2t)^2a_n^2}{2n^2}\notag\\
&\quad-\lambda_i^2\frac{(2t)^3a_n^3}{3!n^3}+o\left(\frac{\lambda_i^2a_n^4}{n^4}\right)+\lambda_i^3\frac{(2t)^3a_n^3}{3n^3}+o\left(\frac{\lambda_i^3a_n^4}{n^4}\right)\label{exp}\\
&=\lambda_i\frac{2ta_n}{n}+\lambda_i\frac{(2t)^2a_n^2}{2n^2}-\lambda_i^2\frac{(2t)^2a_n^2}{2n^2}+o\left(\lambda_i(1-\lambda_i)\frac{a_n^3}{n^3}\right)\label{upperbound}\\
&=\lambda_i\frac{ta_n}{n}+\lambda_i(1-\lambda_i)\frac{t^2a_n^2}{2n^2}+o\left(\frac{\lambda_i(1-\lambda_i)}{n^3}\right).\notag
\end{align}
Here we used Taylor series for the exponential function at $0$ in \eqref{exp} and the following estimate in equation \eqref{upperbound}.
\begin{align*}
 \lambda_i\frac{(2t)^3a_n^3}{3!n^3}-\lambda_i^2\frac{(2t)^3a_n^3}{3!n^3}+\lambda_i^3\frac{(2t)^3a_n^3}{3n^3}&=\lambda_i^2\frac{(2t)^3a_n^3}{3!n^3}\left(1-3\lambda_i+2\lambda_i^2\right)\\
&=\frac{(2t)^3a_n^3}{3!n^3}\lambda_i(1-\lambda_i)(1-2\lambda_i)\\
&\leq\frac{(2t)^3a_n^3}{3!n^3}\lambda_i(1-\lambda_i).
\end{align*}
We therefore finally get
\begin{align*}
 \log\left(1+\lambda_i\left(e^{2ta_n/n-1}\right)\right)&=\frac{(2t)^2a_n^2\sum_{i=1}^n\lambda_i(1-\lambda_i)}{2n^2}+o\left(\frac{a_n^3\sum_{i=1}^n\lambda_i(1-\lambda_i)}{n^3}\right)\\
&=\frac{t^2a_n^2}{2n}\cdot\frac{\sum_{i=1}^n4f(T^ix)(1-f(T^ix))}{n}+o\left(\frac{a_n^3}{n^2}\cdot\frac{\sum_{i=1}^nf(T^ix)(1-f(T^ix))}{n}\right)
\end{align*}
for all $t\in\mathbb{R}$ and $n\rightarrow\infty$. Thus
\begin{equation}
 \lim_{n\rightarrow\infty}\frac{n}{a_n^2}\log\left(1+\lambda_i\left(e^{2ta_n/n-1}\right)\right)=\frac{t^2a}{2}\label{final}
\end{equation}
since $f\in\mathcal{C}_1(\mathcal{S})$ and $a_n/n\rightarrow 0$ as $n\rightarrow\infty$ by assumption. An application of the G\"{a}rtner-Ellis-Theorem (see Th. 3.2.6 in \cite{dembo}) now yields an MDP for $\frac{1}{a_n}\hat{S}_n$ with speed $\frac{a_n^2}{n}$ and rate function
$$I(t)=\sup_x\left\{xt-\frac{t^2a}{2}\right\}=\frac{t^2}{2a}.$$
\par\hspace*{\fill}$\Box$\par\vskip2ex
\noindent
\textit{Proof of Theorem \ref{th1}:}
We would like to prove a MDP for the sequence of probability measures
\begin{align*}
&\quad\left\{P_{n,x}\left(\frac{M_n-nm}{n^{\alpha}}\in\bullet\,\bigg\vert\,\frac{M_n}{n}\in[m-a,m+a]\right)\right\}_{n\in\mathbb{N}}\\
&=\left\{P_{n,x}\left(\frac{M_n-nm}{n^{\alpha}}\in\bullet\,\bigg\vert\,\frac{M_n-nm}{n^{\alpha}}\in[-an^{1-\alpha},an^{1-\alpha}]\right)\right\}_{n\in\mathbb{N}}
\end{align*}
for some $A=A(m)$ and all $0<a<A$. The Lemmas \ref{lemmaExp1} und \ref{lemmaExp2} state that it suffices to prove a MDP for the sequence of measures 
$$\left\{Q_{n,x}\left(\frac{M_n-nm}{n^{\alpha}}+\frac{W}{n^{\alpha-\frac 12}}\in\bullet\,\bigg\vert\,\frac{M_n-nm}{n^{\alpha}}+\frac{W}{n^{\alpha-\frac 12}}\in[-an^{1-\alpha},an^{1-\alpha}]\right)\right\}_{n\in\mathbb{N}}.$$ 
Lemma \ref{trafo} yields for every Borel set $B$
\begin{align*}
&\quad Q_{n,x}\left(\frac{M_n-nm}{n^{\alpha}}+\frac{W}{n^{\alpha-\frac 12}}\in B\,\bigg\vert\,\frac{M_n-nm}{n^{\alpha}}+\frac{W}{n^{\alpha-\frac 12}}\in[-an^{1-\alpha},an^{1-\alpha}]\right)\\
&=\frac{\int_B\exp(-nG_n(sn^{-(1-\alpha)}))~ds}{\int_{-an^{1-\alpha}}^{an^{1-\alpha}}\exp(-nG_n(sn^{-(1-\alpha)}))~ds}.
\end{align*}
We will prove a MDP for this sequence of measures via \textit{Laplace principle}. Theorem 1.2.3 in \cite{dupuis} states that it satisfies a MDP with the respective speed and rate funtion if and only if it satisfies the \textit{Laplace principle}. So let $\Psi\in\mathcal{C}_b(\mathbb{R})$ be a continuous and bounded function. To varify the Laplace principle we have to show that
\begin{align}
&\quad\lim_{n\rightarrow\infty}\frac{1}{n^{1-2k+2k\alpha}}\log\int_{\mathbb{R}}\exp\left[n^{1-2k+2k\alpha}\Psi(s)\right]~Q_{n,x}\left(Y_n\in ds\,\vert\,Y_n\in[-an^{1-\alpha},an^{1-\alpha}]\right)\notag\\
&=\sup_{s\in\mathbb{R}}\left\{\Psi(s)-\lambda\frac{s^{2k}}{(2k)!}\right\}\label{dasZiel},
\end{align}
where we used the abbreviation $Y_n:=\frac{M_n-nm}{n^{\alpha}}+\frac{W}{n^{\alpha-\frac 12}}$. We substitute the density of $Q_{n,x}(Y_n\in\bullet\,\vert\,Y_n\in[-an^{1-\alpha},an^{1-\alpha}])$ on the left hand side of (\ref{dasZiel}) and thus have to analyze the following object:
$$\frac{1}{n^{1-2k+2k\alpha}}\log\left\{\frac{\int_{-an^{1-\alpha}}^{an^{1-\alpha}}\exp\left[n^{1-2k+2k\alpha}\Psi(s)-nG_n(m+sn^{-(1-\alpha)})\right]~ds}{\int_{-an^{1-\alpha}}^{an^{1-\alpha}}\exp\left[-nG_n(m+sn^{-(1-\alpha)})\right]~ds}\right\},$$
or equivalently
\begin{equation}
\label{zielausdruckLaplace}
\frac{1}{n^{1-2k+2k\alpha}}\log\left\{\frac{\int_{-an^{1-\alpha}}^{an^{1-\alpha}}\exp\left[n^{1-2k+2k\alpha}\Psi(s)-n(G_n(m+sn^{-(1-\alpha)})-G_n(m))\right]~ds}{\int_{-an^{1-\alpha}}^{an^{1-\alpha}}\exp\left[-n(G_n(m+sn^{-(1-\alpha)})-G_n(m))\right]~ds}\right\}.
\end{equation}
We consider the nominator and the denominator in \eqref{zielausdruckLaplace} seperately.\par\noindent
Lemma \ref{hauptlemma} states that there exist $r>0$, $N\in\mathbb{N}$ and a polynomial $H(s)=\frac{\lambda}{2}\frac{s^{2k}}{(2k)!}-\sum_{j=1}^{2k}\left|s\right|^j$ such that for all $n\geq N$ and for all $s$ with $\left|s\right|<rn^{1-\alpha}$ the following estimate holds:
\begin{equation}
 \frac{1}{n^{2k(\alpha-1)}}(G_n(m+sn^{-(1-\alpha)})-G_n(m))\geq H(s).\label{item1}
\end{equation}
We choose $A(m):=r$. Since the leading coefficient of $H$ is positive, $H(s)\rightarrow\infty$ for $\left|s\right|\rightarrow\infty$. Since $H(s)\rightarrow\infty$ and $\lambda\frac{s^{2k}}{(2k)!}\rightarrow\infty$ for $\left|s\right|\rightarrow\infty$, there exists $M>0$ such that $\sup_{\left|s\right|>m}\{\Psi(s)-H(s)\}\leq -\left|\Delta\right|-1$
and the supremum of $\Psi-\lambda\frac{s^{2k}}{(2k)!}$ over $\mathbb{R}$ is attained on the interval $[-M,M]$. This and item \eqref{item1} together imply that for all $0<a<A$ and for all $n\geq N$ with $an^{(1-\alpha)}>M$ it holds
\begin{align}
&\quad\sup_{\{M<\left|s\right|<an^{(1-\alpha)}\}}\left\{n^{1-2k+2k\alpha} \Psi(s)- n(G_n(m+sn^{-(1-\alpha)})-G_n(m))\right\}\notag\\
&=\sup_{\{M<\left|s\right|<an^{(1-\alpha)}\}}\left\{n^{1-2k+2k\alpha} \left[\Psi(s)- \frac{1}{n^{2k(\alpha-1)}}(G_n(m+sn^{-(1-\alpha)})-G_n(m))\right]\right\}\notag\\
&\leq-n^{1-2k+2k\alpha}\left(\left|\sup_{s\in\mathbb{R}}\{\Psi(s)-\lambda\frac{s^{2k}}{(2k)!}\}\right|+1\right).\label{item2}
\end{align}
Now let $0<a<A$, with $A=A(m)$ being the constant chosen above. Lemma \ref{hauptlemma} then implies that for all $\delta>0$ and $n$ large enough
$$\left|\frac{1}{n^{2k(\alpha-1)}}(G_n(m+sn^{-(1-\alpha)})-G_n(m))-\lambda\frac{s^{2k}}{(2k)!}\right|<\delta,$$
for all $s\in[-M,M]$, where $M$ is the constant chosen in \eqref{item2}. Thus
\begin{align}
&\quad\exp\left(-n^{1-2k+2k\alpha}\cdot \delta\right)\int_{\{\left|s\right|\leq M\}}\exp\left(n^{1-2k+2k\alpha}\left[\Psi(s)-\lambda\frac{s^{2k}}{(2k)!}\right]\right)~ds\notag\\
&\leq\int_{\{\left|s\right|\leq M\}}\exp\left(n^{1-2k+2k\alpha}\left[\Psi(s)-\frac{1}{n^{2k(\alpha-1)}}(G_n(m+sn^{-(1-\alpha)})-G_n(m))\right]\right)~ds\notag\\
&\leq\exp\left(n^{1-2k+2k\alpha}\cdot\delta\right)\int_{\{\left|s\right|\leq M\}}\exp\left(n^{1-2k+2k\alpha}\left[\Psi(s)-\lambda\frac{s^{2k}}{(2k)!}\right]\right)~ds\label{zentralAbsch1}
\end{align}
Estimate \eqref{item2} implies for $n$ large enough
\begin{align}
&\quad\int_{\{M<\left|s\right|<an^{1-\alpha}\}}\exp\left(n^{1-2k+2k\alpha} \Psi(s)- n(G_n(m+sn^{-(1-\alpha)})-G_n(m))\right)~ds\notag\\
&\leq 2an^{1-\alpha}\exp\left(-n^{1-2k+2k\alpha}\left(\left|\sup_{s\in\mathbb{R}}\{\Psi(s)-\lambda\frac{s^{2k}}{(2k)!}\}\right|+1\right)\right).\label{zentralAbsch2}
\end{align}
For sufficiently large $n$ we get via the estimates \eqref{zentralAbsch1} and \eqref{zentralAbsch2} 
\begin{align*}
&\quad\exp\left(-n^{1-2k+2k\alpha}\cdot \delta\right)\int_{\{\left|s\right|\leq M\}}\exp\left(n^{1-2k+2k\alpha}\left[\Psi(s)-\lambda\frac{s^{2k}}{(2k)!}\right]\right)~ds\notag\\
&\leq\int_{-an^{1-\alpha}}^{an^{1-\alpha}}\exp\left[n^{1-2k+2k\alpha}\Psi(s)-n(G_n(m+sn^{(1-\alpha)})-G_n(m))\right]~ds\\
&\leq\exp\left(n^{1-2k+2k\alpha}\cdot\delta\right)\int_{\left|s\right|\leq M}\exp\left(n^{1-2k+2k\alpha}\left[\Psi(s)-\lambda\frac{s^{2k}}{(2k)!}\right]\right)~ds\\
&+2an^{1-\alpha}\exp\left(-n^{1-2k+2k\alpha}\left(\left|\sup_{s\in\mathbb{R}}\{\Psi(s)-\lambda\frac{s^{2k}}{(2k)!}\}\right|+1\right)\right).
\end{align*}
Lemma \ref{laplaceLemma} and the fact that the supremum of $\Psi-\lambda\frac{s^{2k}}{(2k)!}$ is attained on the interval $[-M,M]$ now yield
\begin{align*}
&\quad\lim_{n\rightarrow\infty}\frac{1}{n^{1-2k+2k\alpha}}\log\int_{\{\left|s\right|\leq M\}}\exp\left(n^{1-2k+2k\alpha}\left[\Psi(s)-\lambda\frac{s^{2k}}{(2k)!}\right]\right)~ds\\
&=\sup_{\{\left|s\right|\leq M\}}\left\{\Psi(s)-\lambda\frac{s^{2k}}{(2k)!}\right\}\\
&=\sup_{s\in\mathbb{R}}\left\{\Psi(s)-\lambda\frac{s^{2k}}{(2k)!}\right\}.
\end{align*}
Therefore
\begin{align*}
&\quad\sup_{s\in\mathbb{R}}\left\{\Psi(s)-\lambda\frac{s^{2k}}{(2k)!}\right\}-\delta\\
&\leq\liminf_{n\rightarrow\infty}\frac{1}{n^{1-2k+2k\alpha}}\log\int_{-an^{1-\alpha}}^{an^{1-\alpha}}\exp\left[n^{1-2k+2k\alpha}\Psi(s)-n(G_n(m+sn^{-(1-\alpha)})-G_n(m))\right]~ds\\
&\leq\limsup_{n\rightarrow\infty}\frac{1}{n^{1-2k+2k\alpha}}\log\int_{-an^{1-\alpha}}^{an^{1-\alpha}}\exp\left[n^{1-2k+2k\alpha}\Psi(s)-n(G_n(m+sn^{-(1-\alpha)})-G_n(m))\right]~ds\\
&\leq\sup_{s\in\mathbb{R}}\left\{\Psi(s)-\lambda\frac{s^{2k}}{(2k)!}\right\}+\delta.
\end{align*}
Since this equality is valid for all $\delta>0$ we finally get
\begin{align*}
&\quad\lim_{n\rightarrow\infty}\frac{1}{n^{1-2k+2k\alpha}}\log\int_{-an^{1-\alpha}}^{an^{1-\alpha}}\exp\left[n^{1-2k+2k\alpha}\Psi(s)-n(G_n(m+sn^{-(1-\alpha)})-G_n(m))\right]~ds\\
&=\sup_{s\in\mathbb{R}}\left\{\Psi(s)-\lambda\frac{s^{2k}}{(2k)!}\right\}.
\end{align*}
Considering the special case $\Psi\equiv 0$ yields
\begin{align*}
&\quad\lim_{n\rightarrow\infty}\frac{1}{n^{1-2k+2k\alpha}}\log\int_{-an^{1-\alpha}}^{an^{1-\alpha}}\exp\left[-n(G_n(m+sn^{-(1-\alpha)})-G_n(m))\right]~ds\\
&=\sup_{s\in\mathbb{R}}\left\{-\lambda\frac{s^{2k}}{(2k)!}\right\}\\
&=0.
\end{align*}
Taking these two limits together we get \eqref{dasZiel}. Now the Lemmas \ref{lemmaExp1} and \ref{lemmaExp2} can be applied in the cases $k\geq 2$ and $k=1$ respectively. This yields the assertion.\par\hspace*{\fill}$\Box$\par\vskip2ex
\noindent
\textit{Proof of Theorem \ref{th1a}:} The Lemmas \ref{lemmaExp1} and \ref{lemmaExp2} state that it suffices to prove a MDP for the seqence of measures
$$\left\{Q_{n,x}\left(\frac{M_n-nm}{n^{\alpha}}+\frac{W}{n^{\alpha-\frac 12}}\in\bullet\,\right)\right\}_{n\in\mathbb{N}}.$$ 
By Theorem \ref{th1} we already know that for some $A(m)$ and all $0<a<A(m)$ the sequence
$$\left\{Q_{n,x}\left(\frac{M_n-nm}{n^{\alpha}}+\frac{W}{n^{\alpha-\frac 12}}\in\,\bullet\;\bigg\vert\,\frac{M_n-nm}{n^{\alpha}}+\frac{W}{n^{\alpha-\frac 12}}\in[-an^{1-\alpha},an^{1-\alpha}]\right)\right\}_{n\in\mathbb{N}}$$
obeys a MDP with the above speed and rate function. We prove the MDP again via Laplace methode, i.e. we consider for $\Psi\in\mathcal{C}_b(\mathbb{R})$
\begin{align*}
&\quad\frac{1}{n^{1-2k+2k\alpha}}\log\int_{\mathbb{R}}\exp\left[n^{1-2k+2k\alpha}\Psi(s)\right]~dQ_{n,\beta}(s)\\
&=\frac{1}{n^{1-2k+2k\alpha}}\log\left\{\frac{\int_{\mathbb{R}}\exp\left[n^{1-2k+2k\alpha}\Psi(s)-nG_n(m+sn^{-(1-\alpha)})\right]~ds}{\int_{\mathbb{R}}\exp\left[-nG_n(m+sn^{-(1-\alpha)})\right]~ds}\right\}
\end{align*}
or equivalently
$$\frac{1}{n^{1-2k+2k\alpha}}\log\left\{\frac{\int_{\mathbb{R}}\exp\left[n^{1-2k+2k\alpha}\Psi(s)-n(G_n(m+sn^{-(1-\alpha)})-G_n(m))\right]~ds}{\int_{\mathbb{R}}\exp\left[-n(G_n(m+sn^{-(1-\alpha)})-G_n(m))\right]~ds}\right\}.$$
und we study the nominator and the denominator separately. To this end we apply the conditioned version of the MDP as stated in Theorem \ref{th1} and the remaining work consists in controlling the missing integrals. For that purpose we make use of the assumed uniqueness of the global minimum of $G$. Due to the lack of space we write
$$(\dots):=\exp\left[n^{1-2k+2k\alpha}\Psi(s)-n(G_n(m+sn^{-(1-\alpha)})-G_n(m))\right].$$
We then get
\begin{align}
&\quad\frac{1}{n^{1-2k+2k\alpha}}\log\int_{\mathbb{R}}(\dots)~ds\notag\\
&=\frac{1}{n^{1-2k+2k\alpha}}\log\bigg\{\int_{[-an^{1-\alpha},an^{1-\alpha}]}(\dots)~ds+\int_{\left|s\right|>an^{1-\alpha}}(\dots)~ds\bigg\},
\intertext{thus}
&\quad\liminf_{n\rightarrow\infty}\frac{1}{n^{1-2k+2k\alpha}}\log\int_{[-an^{1-\alpha},an^{1-\alpha}]}(\dots)~ds\notag\\
&\leq\limsup_{n\rightarrow\infty}\frac{1}{n^{1-2k+2k\alpha}}\log\bigg\{\int_{[-an^{1-\alpha},an^{1-\alpha}]}(\dots)~ds+\int_{\{\left|s\right|>an^{1-\alpha}\}}(\dots)~ds\bigg\}\notag\\
&=\max\bigg\{\limsup_{n\rightarrow\infty}\frac{1}{n^{1-2k+2k\alpha}}\log\int_{[-an^{1-\alpha},an^{1-\alpha}]}(\dots)~ds,\notag\\
&\quad\limsup_{n\rightarrow\infty}\frac{1}{n^{1-2k+2k\alpha}}\log\int_{\left\{\left|s\right|>an^{1-\alpha}\right\}}(\dots)~ds\bigg\}.\label{einkaestelung}
\end{align}
Theorem \ref{th1} implies 
$$\quad\lim_{n\rightarrow\infty}\frac{1}{n^{1-2k+2k\alpha}}\log\int_{[-an^{1-\alpha},an^{1-\alpha}]}(\dots)~ds=\sup_{s\in\mathbb{R}}\left\{\Psi(s)-\lambda\frac{s^{2k}}{(2k)!}\right\}.$$
We now show that the second argument of the maximum converges to $-\infty$ as $n$ goes to $\infty$. Thereby we make use of the uniqueness of the global minimum $m$ of $G$. Lemma \ref{dreiLemma} \textit{(iii)} implies
\begin{align}
&\quad\frac{1}{n^{1-2k+2k\alpha}}\log\int_{\left\{\left|s\right|>an^{1-\alpha}\right\}}(\dots)~ds\notag\\
&=\frac{1}{n^{1-2k+2k\alpha}}\log\int_{\left\{\left|s\right|>an^{1-\alpha}\right\}}\exp\left[n^{1-2k+2k\alpha}\Psi(s)-n(G_n(m+sn^{-(1-\alpha)})-G_n(m))\right]~ds\notag\\
&=\frac{1}{n^{1-2k+2k\alpha}}\log\int_{\left\{\left|t-m\right|>a\right\}}\exp\left[n^{1-2k+2k\alpha}\Psi((t-m)n^{1-\alpha})-n(G_n(t)-G_n(m))\right]n^{1-\alpha}~dt\notag\\
&\leq\frac{1}{n^{1-2k+2k\alpha}}\log\left\{\exp\left[n^{1-2k+2k\alpha}\|\Psi\|_{\infty}\right]\cdot n^{1-\alpha}e^{nG_n(m)}\underbrace{\int_{\left\{\left|t-m\right|>a\right\}}\exp(-nG_n(t))~dt}_{=\mathcal{O}(e^{-n\epsilon-n g})}\right\}\notag\\
&=\|\Psi\|_{\infty}+\frac{1}{n^{1-2k+2k\alpha}}\log\left\{n^{1-\alpha}\mathcal{O}\left(\exp\left[-n(\epsilon+(g-G_n(m)))\right]\right)\right\}\notag\\
&=\|\Psi\|_{\infty}+\underbrace{\frac{\log(n^{1-\alpha})}{n^{1-2k+2k\alpha}}}_{\rightarrow 0}+\frac{1}{n^{1-2k+2k\alpha}}\log\left\{\mathcal{O}\left(\exp\left[-n(\epsilon+\underbrace{(g-G_n(m))}_{\rightarrow 0})\right]\right)\right\}\notag\\
&\underset{n\rightarrow\infty}{\longrightarrow}-\infty,\label{vernachl}
\end{align}
since the set $\left\{\left|t-m\right|>a\right\}$ does not contain a minimum of $G$. This together with inequality \eqref{einkaestelung} yields 
\begin{align*}
\quad\lim_{n\rightarrow\infty}\frac{1}{n^{1-2k+2k\alpha}}\log\int_{\mathbb{R}}(\dots)~ds&=\lim_{n\rightarrow\infty}\frac{1}{n^{1-2k+2k\alpha}}\log\int_{[-an^{1-\alpha},an^{1-\alpha}]}(\dots)~ds\\
&=\sup_{s\in\mathbb{R}}\left\{\Psi(s)-\lambda\frac{s^{2k}}{(2k)!}\right\}.
\end{align*}
Considering the special case $\Psi\equiv 0$ implies
\begin{align*}
&\quad\lim_{n\rightarrow\infty}\frac{1}{n^{1-2k+2k\alpha}}\log\int_{\mathbb{R}}\exp\left[-n(G_n(m+sn^{-(1-\alpha)})-G_n(m))\right]~ds\\
&=\sup_{s\in\mathbb{R}}\left\{-\lambda\frac{s^{2k}}{(2k)!}\right\}=0.
\end{align*}
These two limits together yield \eqref{dasZiel}. By Lemma \ref{lemmaExp1} and Lemma \ref{lemmaExp2} we get the assertion for the case $k\geq 2$ and $k=1$ respectively.\par\hspace*{\fill}$\Box$\par\vskip2ex\par
For illustrational purposes we give a concrete example of a dynamical system and compute the rate function of the respective MDP explicitly. Note that this dynamical system has already been considered in \cite{guillotin}. We refer to this paper in order to check that our conditions imposed in Theorem \ref{th1a} hold true for this concrete example. 
\section{Example: irrational rotation on the torus}
In \cite{guillotin} the authors consider the dynamical system $(\mathbb{T}, \mathcal{B}(\mathbb{T}), \lambda_{\mathbb{T}}, T_{\alpha})$. There $\mathbb{T}=\mathbb{R}/\mathbb{Z}=[0,1[$ denotes the one-dimensional torus, $\lambda_{\mathbb{T}}$ the restricted Lebesgue-measure and $T_{\alpha}$ the irrational rotation with angle $\alpha$ of type $\eta$ (see Definition 5.3 in \cite{guillotin}), i.e. $x\mapsto x+\alpha \text{ mod } 1$.\\
Let $f(x)=x$ be the identity on $\mathbb{T}$. In \cite{guillotin} it is proved that the Curie-Weiss model with dynamical external field according to this dynamical system exhibits a phase transition at the critical temperature $\beta_c=\frac{3}{2J}$.
\begin{Theorem}[MDP for the irrational rotation on the one-dimensional torus]
For the above defined dynamical system the following assertions hold true.
\begin{enumerate}
\item For $\beta<\beta_c$, $\eta<2$ and for all $\frac{1}{2}<\alpha<1$ the sequence of measures 
$$\left\{P_{n,x}\left(\frac{M_n-nm}{n^{\alpha}}\in\bullet\right)\right\}_{n\in\mathbb{N}}$$
satisfies a MDP with scale $n^{2\alpha-1}$ and rate function
$$I(z)=\frac{z^2}{2\sigma^2},$$
where $\sigma^2=\frac{2}{3-2\beta J}$.
\item For $\beta=\beta_c$, $\eta<4/3$ and for all $\frac{1}{2k}<\alpha<1$ the sequence of measures $$\left\{P_{n,x}\left(\frac{M_n-nm}{n^{\alpha}}\in\bullet\right)\right\}_{n\in\mathbb{N}}$$
satisfies a MDP with scale $n^{1-2k+2k\alpha}$ and rate function
$$I(z)=\frac{9}{80}z^{4}.$$
\end{enumerate}
\end{Theorem}


\end{document}